\documentclass[review]{elsarticle}
\usepackage{hyperref}
\usepackage{amsfonts}

\begin{document}

\newtheorem{thm}{Theorem}
\newtheorem{lem}[thm]{Lemma}
\newtheorem{cor}[thm]{Corollary}
\newtheorem{pro}[thm]{Proposition}
\newtheorem{conj}[thm]{Conjecture}
\newdefinition{dfn}[thm]{Definition}
\newproof{prf}{Proof}

%\begin{frontmatter}

\title{On Structural Descriptions of Lower Ideals of Series Parallel Posets}
\author[cja]{Christian Joseph Altomare}
\ead{caltomare@towson.edu}
\address{The Ohio State University, 231, West 18th Avenue,
Columbus, Ohio, United States}

%undone
\begin{abstract}
In this paper we give an algorithm to determine, for any given suborder
closed class of series-parallel posets, a structure theorem for the class.
We refer to these structure theorems as structural descriptions.
\end{abstract}

\maketitle

%\end{frontmatter}

%undone

\section{Introduction}

\def\rep{representation}
\def\bs{\hbox{BS}(P)}
\def\bss{\hbox{BS}}
\def\t{\hbox{Top}}
\def\b{\hbox{Bottom}}
\def\p{\prec}
\def\cl{\hbox{cl}}
\def\forb{\hbox{Forb}}
\def\cs{P_1\p\cdots\p P_n}
\def\as{P_1\oplus\cdots\oplus P_n}
\def\lexpi{\bigoplus_{\le_I}P_i}
\def\bit{(I,\le_I,f)}
\def\xs{X_1\p\cdots\p X_n}
\def\a{\oplus}
\def\bo{\bigoplus}
\def\o{\oplus}
\def\fg{\forb(\G;P_1,\ldots,P_k)}
\def\ab{(A,B)}
\def\lcis{left cell ideal set}
\def\rcis{right cell ideal set}
\def\ln{\hbox{lcis}}
\def\rn{\hbox{rcis}}
\def\bsspg{\bss(\G;P_1,\ldots,P_k)}
\def\lcl{\hbox{lcl}(g;P_1,\ldots,P_k)}
\def\rcl{\hbox{rcl}(g;P_1,\ldots,P_k)}
\def\potk{P_1,\ldots,P_k}
\def\potn{P_1,\ldots,P_n}
\def\s{\hbox{spl}}
\def\G{\Gamma}

\subsection{Background}

%undone def of finite struct thm is off
%current def would call it finite even if it's infinite at a more deeply
%nested level

%undone possible confusion as you say all partial orders in this paper
%are finite, yet the wqo of SP orders is infinite. This particular wqo
%is also a partial order, so it could confuse a reader.

%undone don't just state as nigussie's.
%cite nigussie-robertson AND robertson-seymour-thomas

%undonee combine defs of X_1<...X_n and P_1<...<P_n

%undone skip MUCH of technical lemma section by finding appropriate references

Many important theorems in combinatorics characterize a class by forbidden
subobjects of some kind. This is a description of the class ``from the
outside'', by what is not inside it. An example is Wagner's reformulation
\cite{wagner}
of Kuratowski's Theorem \cite{kuratowski}
stating that a graph is planar iff it has no $K_5$
minor and no $K_{3,3}$ minor. To be a good characterization, the list of
forbidden objects should be finite. Well quasi order theorems such as the
Graph Minor Theorem \cite{graph_minor_theorem}
state that for certain classes of objects,
there is always such a finite description ``from the outside''.

Just as important are those theorems that
characterize a class ``from the inside'' by giving some set of starting
objects and some set of construction rules. As a simple example, consider
(graph theoretic) trees. Each tree is either a single point graph or may
be obtained from two smaller, disjoint trees by adding an edge between
the trees. Therefore a simple structure theorem for this class would have
the single point graph as the only starting graph and joining two disjoint
graphs by an edge as the sole construction rule.

To be a good characterization,
we again hope that it is in some sense finite. First, there should be only
finitely many construction rules. We can not necessarily demand there are
only finitely many starting objects. We may however demand that at least we
start with only finitely many families, and that each such family has some
sort of finite description as well.

Analogous to the Graph Minor Theorem
and other well quasi order theorems stating that in many cases
there is always a finite description from the outside, it was asked
if it could be shown in an equally general setting that there is always
a finite description from the inside with finitely many starting families,
each itself finitely described, and finitely many construction rules.

As it turns out, this appears to be far more difficult. This line of research
was first pursued by Robertson, Seymour, and Thomas in \cite{rst} for trees
under the topological minor relation.
In \cite{nigussie_robertson}, Nigussie and Robertson build on \cite{rst} and
correct some technical errors contained therein. In \cite{nigussie},
Nigussie gives an algorithm that finds a structure theorem for an arbitrary
topological minor closed property of trees. Nigussie's algorithm is efficient
enough in practice that structure theorems can be computed by hand with pen
and paper that are not at all obvious without the algorithm.
We follow the convention of referring to these structure theorems
as structural descriptions. The distinction
we make is that we use the term structure theorem  informally, while we see
 structural description as a technical term defined in \cite{nigussie} for
trees under topological minor and below for series-parallel orders under
suborder.

Attempts have been made by various researchers to generalize these
results to other classes of graphs, in particular series-parallel graphs.
Thus far, no such attempt has succeeded. While many specific graph
structure theorems are known, the tree result is to date the only one
that allows the automatic computation of a structure theorem for any graph
property in a nontrivial, infinite class of properties.

It is key that rooted trees are used in \cite{rst}, \cite{nigussie_robertson},
and \cite{nigussie}. Rooted trees are as much partial orders as they are
graphs, and we view Nigussie's algorithm not just as a graph algorithm, but as
a partial order algorithm. It is thus natural to ask for algorithms similar to
Nigussie's for classes of partial orders larger than the class of trees.
In this paper, we prove an analogous result for series-parallel partial
orders by giving a finite structural description for each suborder closed class
of series-parallel orders. More precisely, we give an algorithm that takes
as input a suborder closed class of series-parallel
orders described by forbidden suborders, and
which gives as output a finite structural description for that class.

In our context, a structural description will turn out to be a finite set of
labeled partial orders. The labels will be families already constructed.
Each labeled partial order in the structural description for a class will
represent one family or construction rule. Roughly, the labels tell what
we are allowed to put in and the partial orders themselves tell us how
we are allowed to piece together what we do put in.

\section{Basic Definitions and Conventions}

A partial order is a (possibly empty)
set $P$ together with a reflexive, antisymmetric,
transitive binary relation $\le$ on $P$. All partial orders in this paper
are assumed to be finite. (The only exception to this is that {\it classes
of} partial orders we consider are usually infinite, and this class together
with the suborder relation is in fact a partial order. This exception causes
no confusion as it is clear in each case whether we are dealing with a
partial order or an infinite family of them.)
Points $x,y$ in a partial order
$(P,\le)$ are comparable if $x\le y$ or $y\le x$. Otherwise $x$ and $y$
are called incomparable, which we write as $x|y$. A chain is a partial order
such that any two points are comparable.
An antichain is a partial order such that any two points are incomparable.

A lower ideal of partial orders is a family of partial orders that
is closed under taking suborders.
Given partial orders $P$ and $Q$, we say that $P$ is $Q$-free if $P$ has
no suborder isomorphic to $Q$. Given a set $F$ of partial orders, we say
that $P$ is $F$-free if $P$ is $Q$-free for each $Q$ in $F$. A lower ideal $L$
is said to be $Q$-free or $F$-free if each partial order in $L$ is
$Q$-free or $F$-free, respectively.
A forbidden suborder of a lower ideal $L$ is a suborder minimal partial order
$P$ such that $L$ is $P$-free.

The papers \cite{rst}, \cite{nigussie_robertson}, and \cite{nigussie}
use tree sums to construct new trees from old.
For our purposes, tree sums are not sufficient. The correct
generalization to our context is partial order lexicographic sums. We call
partial orders $(P_i,\le_i)$ and $(P_j,\le_j)$ disjoint if $P_i$ and
$P_j$ are disjoint.

\begin{dfn}
Let $\{(P_i,\le_i)\}_{i\in I}$ be a family of pairwise disjoint partial
orders and let $(I,\le_I)$ be a partial order on $I$. Then the lexicographic
sum $\bigoplus_{\le_I}P_i$ is defined as the unique partial order
$(\bigcup_{i\in I}P_i,\le)$ such that the following conditions hold:
\begin{enumerate}
\item Given $i$ in $I$ and $x$ and $y$ in $P_i$, we have
$x\le y$ iff $x\le_i y$.
\item Given distinct $i, j$ in $I$, if $i\le_I j$, then $x\le y$ for all
$x$ in $P_i$ and $y$ in $P_j$.
\item Given distinct $i, j$ in $I$, if $i$ and $j$ are $\le_I$ incomparable,
then $x$ and $y$ are $\le$ incomparable for all
$x$ in $P_i$ and $y$ in $P_j$.
\end{enumerate}
\end{dfn}

It is a simple exercise to show that the above three conditions indeed
uniquely determine a partial order on $\bigcup_{i\in I}P_i$.
We call $(I,\le_I)$ the outer partial order of the lexicographic sum.
Each $P_i$ is called the inner partial order corresponding to $i$. The
lexicographic sum is therefore a partial order on the union of the inner
partial orders.
We call the partition of $\bigoplus_{\le_I}P_i$ into
the inner partial orders $P_i$ a lexicographic partition. It is simple
to show that a partition of a partial order is lexicographic iff for
any two distinct cells $C_1$ and $C_2$ of the partition, either all
elements of $C_1$ precede all elements of $C_2$,
all elements of $C_2$ precede all elements of $C_1$,
or all elements of $C_1$ and $C_2$ are incomparable.
In this case, the outer partial order is uniquely determined in the
obvious way.

We call a lexicographic partition nontrivial if there are at least two
cells and each cell is nonempty.
We call a lexicographic partition a chain partition if the corresponding
outer partial order is a chain. Similarly for antichain partitions.
We call a lexicographic sum a chain sum or antichain sum if the
corresponding partition is nontrivial and the
outer partial order is a chain or antichain, respectively. We denote
by $P_1\p\cdots\p P_n$ the chain sum of partial orders $P_1,\ldots,P_n$
such that for $1\le i<j\le n$, every $x$ in $P_i$ is less than every
$y$ in $P_j$. We denote
by $P_1\oplus\cdots\oplus P_n$ the antichain sum of partial orders
$P_1,\ldots,P_n$
such that for all $i\not=j$, every $x$ in $P_i$ is incomparable to every
$y$ in $P_j$.

The comparability graph of a partial order $P$ is the graph whose vertices
are the points of $P$ and such that two points $x$ and $y$ are adjacent
iff they are comparable in $P$.
A component of $P$ is a component of the comparability graph. An
anticomponent is a component of the similarly defined incomparability graph.
If $P$ is a  chain sum, we note that $P$
then has a unique finest chain partition, which is just the partition into
anticomponents.
If $P=P_1\p \cdots\p P_n$ and $\{P_1,\ldots,P_n\}$ is a finest chain partition
with $n\ge 2$,
then we call $P_1\p \cdots\p P_n$ a finest chain representation of $P$.
A similar statement holds for  antichain sums and components, and
we then similarly call $P_1\bigoplus\cdots\bigoplus P_n$ a finest antichain
representation of $P$ for $n\ge 2$.

%undone italicize definitions everywhere
%undone perhaps find basic on SP orders written somewhere and just QUOTE
%THAT REFERENCE, making a couple paragraphs CLAIMING the SP order facts
%needed, so the paper is shorter.

A partial order is a series-parallel partial order, or SP order, if it
is contained in the smallest class of partial orders containing the
empty and single point partial orders and closed under
chain and antichain sums.
We note that for each SP order $P$, exactly one of the following holds:
$P$ is empty, $P$ is a single point, $P$ is a chain sum, or $P$ is an
antichain sum.
We will make use of the simple but important fact that a suborder of an
SP order is also an SP order.
It is also worth noting  that a finite partial order is an SP
order iff it is $N$-free, where $N$ is the
partial order on points $a,b,c,d$ such that $a<b$, $b>c$, $c<d$, and
all others pairs of points are incomparable \cite{poset_unique_lex_decomp},
though we do not make use of this fact.

Since all our ideals in this paper are lower ideals
of SP orders, from now on we simply call these lower ideals.
A proper lower ideal is a lower ideal that is strictly contained
in the set of all SP orders.
A nontrivial lower ideal is one that contains at least one nonempty
partial order.
Our goal in this paper is to give a structural description for an arbitrary
nontrivial, proper lower ideal. More precisely, we give a
recursive procedure that takes as input a nontrivial, proper lower ideal,
which gives as output a structural description for that lower ideal.
This procedure is entirely constructive, and a program could be written
to implement it, though algorithmic questions are not our focus.

A structural description, for us, will turn out to be a finite set of labeled
SP orders. The labels tell us which objects we may use to construct,
and the orders themselves tell us in which ways we may put these together.
We now start to make this intuition more precise.

A labeled partial order is a triple $(I,\le_I,f)$, where
$(I,\le_I)$ is a partial order and $f$ is a function
with domain $I$. We think of $f$ as the labeling function. We sometimes
write $I_f$ for this labeled partial order when $\le_I$ is clear from context.
A bit is a labeled SP order such that each label is a lower ideal
or the symbol R. We call a point $i$ in a bit $(I,\le_I,f)$ an ideal labeled
point if $f(i)$ is a lower ideal. We call $i$ an $R$ labeled point if
$f(i)=R$. A recursive bit is a bit with at least one $R$ labeled point.
A nonrecursive bit is a bit with no $R$ labeled points.
The two point chain with both points labeled $R$ is denoted by $R_C$.
The two point antichain with both points labeled $R$ is denoted by $R_A$.

We now tell how to assign to each set $S$ of bits the lower
ideal $L(S)$ that $S$ is said to generate. Given a set $S$ of bits and a set
$X$ of partial orders, we say that $X$ is $S$-bit closed if $X$ contains all
lexicographic sums of the form
$\bigoplus_{\le_I}P_i$ such that $(I,\le_I,f)$ is a bit in $S$,
the partial order $P_i$ is contained in the lower ideal $f(i)$ for
each ideal labeled point $i$ in $I$,
and $P_i$ is contained in $X$ itself for each $R$ labeled point $i$ in $I$.
The $S$-bit closure of $X$ is the smallest $S$-bit
closed set containing $X$. Given a set $S$ of bits, we define the lower
ideal $L(S)$ generated by $S$ as the $S$-bit closure of the set
containing the empty partial order, the one point partial order, and
no other partial orders.

Given a bit $\bit$, we say that $X$ is $\bit$-bit closed if $X$ is\\
$\{\bit\}$-bit closed.
We will have many occasions to use the following simple lemma,
whose proof is immediate from the definition.

\begin{lem} \label{set_to_bit}
If $S$ is a set of bits and $X$ is a set of partial orders,
then $X$ is $S$-bit closed iff $X$ is $\bit$-bit closed for each
bit $\bit$ in $S$.
\end{lem}

We now define structural descriptions.
We do so by recursively defining structural descriptions
of each nonnegative integer {\it rank}.
The empty set, thought of as an empty set of bits, is the only structural
description of rank $0$.
Assume the structural descriptions of ranks $0,\ldots,n$ are known.
A structural description of rank $n+1$ is a finite set $S$ of finite, labeled
SP orders
such that each label of each bit in $S$ is either the special symbol ``$R$''
or a structural description of rank at most $n$.
A structural description is a structural description of some finite rank.
Note that since we require finiteness at each step, each of our structural
descriptions would be considered a ``finite structural description'' in the informal sense of the term.

A structural description $D$ generates a lower ideal $L(D)$ analogously to
the previous definition for bits. We say this is a structural description
{\it for} $L(D)$ or {\it of } $L(D)$.

Our recursive procedure will take a lower ideal as input and give
a finite structural description as output. We just made precise what form the
output takes. To state the form of the input, we first need
several definitions.
A quasi order is a set $Q$ together with a transitive, reflexive
relation $\le$.
A quasi order is a well quasi order, or WQO, if for
all infinite sequences $q_1,q_2,\ldots$ of points in $Q$, there are
positive integers $i<j$ such that $q_i\le q_j$. A class $\cal{C}$ of partial
orders is then said to be well quasi ordered under suborder if
for each infinite sequence $(P_1,\le_1),(P_2,\le_2),\ldots$, of partial
orders in $\cal{C}$, there are positive integers $i<j$ such that
$(P_i,\le_i)$ is a suborder of $(P_j,\le_j)$.

Given an SP order $P$,
we let $\forb(P)$ be the set of SP orders forbidding $P$ as a suborder.
Given a set $F=\{Q_1,\ldots,Q_k\}$ of SP orders,
we denote the set of SP orders forbidding each $P$ in $F$ as a suborder
by $\forb(F)$ or $\forb(Q_1,\ldots,Q_k)$.
It can be shown that finite SP orders form a WQO under the suborder relation.
Basic WQO theory then implies that for each lower ideal $L$, there is a
finite set $F$ of SP orders such that $L=\forb(F)$
\cite{kruskal_wqo_intro}.
With these facts stated,
we may now express the main result of this paper more precisely; we give
an algorithm that takes a finite set $F$ of SP orders as input and outputs
a structural description $D$ such that $L(D)=\forb(F)$. 

Since our main focus is
combinatorial structure theory, we do not concern ourselves with algorithmic
or complexity theoretic questions. Though such questions may be
interesting, they are simply not our focus here.
We thus present our algorithms in the same
informal style that is common in mathematics.

\section{Technical Lemmas}

We note that the reader familiar with SP orders can likely skim or
skip much of this section. Even readers unfamiliar with SP orders
may find it useful to proceed to the next section and refer back 
to this section as needed.

We call an SP order connected if its comparability graph is connected.
An SP order is anticonnected if its incomparability graph is connected.

\begin{lem}  \label{chain_sum_connected}
Every chain sum is connected. Similarly,
every antichain sum is anticonnected.
\end{lem}

\begin{prf}
Let $P_1\p \cdots \p P_n$ be a chain sum. By definition, we can assume without
loss of generality that $n\ge 2$ and each $P_i$ is nonempty. For $i\not=j$,
each point $x$ of $P_i$ is comparable to each point $y$ in $P_j$ and
hence $x$ and $y$ are adjacent in the
comparability graph. If two points $x$ and $y$ are contained in the same
$P_i$, then choose $i\not=j$ and $z$ in $P_j$. Then $x$ and $y$ are both
adjacent to $z$ and hence in the same component. Therefore given any
points $x$ and $y$ in $P_1\p \cdots\p P_n$, there is a path of length one or
two between $x$ and $y$ in the comparability graph of $P$, and the first
claim of the lemma holds. For the second claim, repeat the same proof
with $P_1\bigoplus\cdots\bigoplus P_n$ and the incomparability graph.
\end{prf}

\begin{lem} \label{comp_in_P_i}
Each component of $P_1\bigoplus\cdots\bigoplus P_n$  is contained in
some $P_i$.
Each anticomponent of $P_1\p \cdots\p  P_n$ is contained in
some $P_i$.
\end{lem}

\begin{prf}
A component of $P_1\bigoplus\cdots\bigoplus P_n$ is connected in the
comparability graph. Since there are no edges from $P_i$ to $P_j$
for $i\not=j$ in the comparability graph, we see that each component
is contained in some $P_i$. The proof of the second claim is analogous.
\end{prf}

\begin{lem} \label{antichain_bit_there}
If $Q$ is a chain sum and $P_i$ is $Q$-free for $i$ in $\{1,\ldots,n\}$,
then $P_1\bigoplus\cdots\bigoplus P_n$ is $Q$-free.
\end{lem}

\begin{prf}
Let $Q$ be a chain sum. It is enough to show that if
$P_1\bigoplus\cdots\bigoplus P_n$ contains $Q$,
then $P_i$ contains $Q$ for some $i$. Since $Q$ is a chain sum,
we know by \ref{chain_sum_connected} that $Q$ is
connected. By \ref{comp_in_P_i}, $Q$ must therefore be contained
in some $P_i$.
\end{prf}

The next lemma is analogous to the previous lemma, and the same proof
goes through mutatis mutandis.

\begin{lem} \label{chain_bit_there}
If $Q$ is an antichain sum and $P_i$ is $Q$-free for $i$ in $\{1,\ldots,n\}$,
then $P_1\p \cdots\p P_n$ is $Q$-free.
\end{lem}

We need several technical lemmas.

\begin{lem} \label{finest_chain_rep_part}
If $\cs$ is a finest chain representation of an SP order $P$,
then each $P_i$ is an antichain sum or a one point partial order.
\end{lem}

\begin{prf}
For each $i$, since $P_i$ is a suborder of an SP order, $P_i$ itself is
an SP order.
Since $\cs$ is a finest chain representation by hypothesis, it follows
by definition of finest chain representation that $P_i$ is not itself
a chain sum, and $P_i$ is therefore a single point or an antichain sum
as claimed.
\end{prf}

The same holds for finest antichain representations. We omit the entirely
analogous proof.

\begin{lem} \label{finest_antichain_rep_part}
If $\as$ is a finest antichain representation of an SP order $P$,
then each $P_i$ is a chain sum or a one point partial order.
\end{lem}

\begin{lem} \label{finest_antichain_rep_Q}
Let $\as$ be a finest antichain representation of a partial order $P$
and let $Q_1\oplus\cdots\oplus Q_k$ be an arbitrary antichain sum.
If $\as$ is a suborder of $Q_1\oplus\cdots\oplus Q_k$,
then for each $i$ with $1\le i\le n$ there is $j$
with $1\le j\le k$ such that $P_i$ is a suborder of $Q_j$.
\end{lem}

\begin{prf}
Choose $i$. Note that $P_i$ is a chain sum or a one point partial order
by \ref{finest_antichain_rep_part}. If $P_i$ is a single point, then $P_i$
is of course contained in some $Q_i$. If $P_i$ is a chain sum, then it is
connected and therefore contained in a component of
$Q_1\oplus\cdots\oplus Q_k$. Since each component of
$Q_1\oplus\cdots\oplus Q_k$ is contained in some $Q_i$, the result
follows.
\end{prf}

The following lemma has a similar proof.

\begin{lem} \label{finest_chain_rep_Q}
Let $\cs$ be a finest chain representation of a partial order $P$
and let $Q_1\p \cdots\p Q_k$ be an arbitrary chain sum. If $\cs$ is a suborder
of $Q_1\p \cdots\p  Q_k$, then for each $i$ with $1\le i\le n$ there is $j$
with $1\le j\le k$ such that $P_i$ is a suborder of $Q_j$.
\end{lem}

\begin{lem} \label{interval_suborder}
Let $\cs$ be a finest chain representation of an SP order $P$
that is contained in the partial order $Q_1\p Q_2$. If the $P_i$
of $\cs$ is contained in $Q_1$ then so is $P_1\p \cdots\p P_i$.
Similarly, 
if the $P_i$ of $\cs$ is contained in $Q_2$ then so is
$P_i\p \cdots\p P_n$.
\end{lem}

\begin{prf}
We prove the first claim. The second is similar. By hypothesis, the
$P_i$ of $\cs$ is a suborder of $Q_1$. Since every point of
$P_1\p \cdots\p P_i$ is less than or equal some point of $P_i$, and since
$Q_1$ is a downward closed subset of $Q_1\p Q_2$ containing $P_i$,
it follows that $P_1\p \cdots\p P_i$ is a suborder of $Q_1$.
\end{prf}

\begin{lem}
If $\cs$ is a finest chain representation that is contained in
the partial order $Q_1\p Q_2$, then one of the following three conditions
holds:
\begin{enumerate}
\item $\cs$ is a suborder of $Q_1$.
\item $\cs$ is a suborder of $Q_2$.
\item There is $i$ with $1\le i<n$ such that
 $P_1\p \cdots\p  P_i$ is a suborder of $Q_1$ and $P_{i+1}\p \cdots\p P_n$ is
 a suborder of $Q_2$.
\end{enumerate}
\end{lem}

\begin{prf}
Since $\cs$ is a finest chain representation
by hypothesis, we know that
each $P_i$ is contained in $Q_1$ or $Q_2$
by \ref{finest_chain_rep_Q}.
If $\cs$ is a suborder of $Q_1$ or $Q_2$ then we are done. Suppose not.
Take the largest $i$ such that $P_i$ is a suborder of $Q_1$. By
\ref{interval_suborder}, we see
that $P_1\p \cdots\p P_i$ is a suborder of $Q_1$. Since $\cs$ is not a suborder
of $Q_1$ by hypothesis, we know that $i<n$. Therefore $P_{i+1}$ is a suborder
of $Q_2$. Again by \ref{interval_suborder},
we see that $P_{i+1}\p \cdots\p P_n$ is a suborder of $Q_2$,
which completes the proof.
\end{prf}

\section{The Main Lemmas}

Given labels $X_1,\ldots,X_n$, we let the notation
$\xs$ denote the $n$ point labeled chain with bottom
point labeled $X_1$, next least point labeled $X_2$, and so on. Note
that $\cs$ defined previously is the chain sum of $n$ partial orders
$P_1,\ldots,P_n$ (which is of course itself a partial order).
On the other hand, $\xs$ denotes a labeled $n$ point chain bit.
As long as the reader keeps this distinction in mind, no confusion arises.
Similarly for the expression $X_1\o\cdots\o X_n$.

\begin{dfn} \label{chain_bit_choice_first_def}
Let $n\ge 2$.
The chain bit set $\bs$ corresponding to a chain $P$ with finest chain \rep\
$\cs$ is defined to be the set of bits $B$ such that one of the following
conditions hold:
\begin{enumerate}
%\item $B=R_A$.
\item $B=R\p\forb(P_n)$.
\item $B=\forb(P_1)\p R$.
\item There is $i$ with $1<i<n$ such that
      $$B=\forb(P_1\p \cdots\p P_i)\p \forb(P_i\p \cdots\p P_n)$$
\end{enumerate}
\end{dfn}

We note that since the finest chain representation is uniquely determined,
the notation $\bs$ is well defined for chain sums $P$.

\begin{lem} \label{one_chain_lem}
Let $n\ge 2$.
If $P$ is an SP order with finest chain \rep\ $\cs$,
then $$\forb(P_1\p \cdots \p P_n)=L(\bs\cup\{R_A\}).$$
\end{lem}

\begin{prf}
Let $S=\bss(P_1\p\cdots\p P_n)\cup\{R_A\}$.
We must show that $\forb(P_1\p\cdots\p P_n)$ is the $S$-bit closure of the
doubleton containing the empty and one point partial orders.
%undone get rid of brak in middle of notation
Since $\forb(P_1\p\cdots\p P_n)$ trivially contains the empty and
one point partial orders, it is enough
to show that $\forb(P_1\p\cdots\p P_n)$ is $S$-bit closed and that every
$S$-bit closed set containing the empty and one point partial orders
has $\forb(P_1\p\cdots\p P_n)$ as a subset.

We first show that $\forb(P_1\p\cdots\p P_n)$ is $S$-bit closed.
By \ref{set_to_bit}, it is enough to show that $\forb(P_1\p\cdots\p P_n)$ is
$\bit$-bit closed for each bit $\bit$ in $S$. We consider four cases.

First, if $\bit$ is $R_A$,
then to show that $\forb(\cs)$ is $\bit$-bit closed is simply to show that
$\forb(\cs)$ is closed under antichain sums. But this is exactly
\ref{antichain_bit_there}.

Second, if $\bit$ is a two point chain with bottom point labeled $R$
and top point labeled $\forb(P_n)$,
then to show that $\forb(\cs)$ is $\bit$-bit closed is to show that
if $Q_1$ is a partial order in $\forb(\cs)$ and $Q_2$ is a partial order
in $\forb(P_n)$, then $Q_1\p Q_2$ forbids $\cs$. Suppose not. Since $Q_1\p Q_2$
contains $\cs$, in particular $Q_1\p Q_2$ contains the top inner part $P_n$
of the chain sum. By \ref{finest_chain_rep_Q},
we see that $P_n$ is a suborder of $Q_1$
or $Q_2$. Since $Q_2$ forbids $P_n$, we know that $P_n$ is a suborder of
$Q_1$. By \ref{interval_suborder},
it follows that $\cs$ is a suborder of $Q_1$,
contrary to hypothesis. This
contradiction shows that $\forb(\cs)$ is $\bit$-bit closed as claimed.

The third case, that $\bit$ is a two point chain with top point labeled $R$
and bottom point labeled $\forb(P_1)$,
is completely analogous to the second case,
and the proof goes through mutatis mutandis.

Fourth, if there is $i$ with $1<i<n$ such that
$\bit$ is a two point chain with bottom point labeled
$\forb(P_1\p \cdots\p P_i)$ and top point labeled $\forb(P_i\p \cdots\p P_n)$,
then to show that $\forb(\cs)$ is $\bit$-bit closed, we must show
that if $Q_1$ is a partial order forbidding
$P_1\p \cdots\p P_i$ and $Q_2$ is a partial order forbidding $P_i\p \cdots\p P_n$,
then $Q_1\p Q_2$ forbids $\cs$. We prove the contrapositive statement,
namely, that if $Q_1\p Q_2$ has a $\cs$ suborder then
$Q_1$ has a $P_1\p \cdots\p P_i$ suborder or
$Q_2$ has a $P_i\p \cdots\p P_n$ suborder. Since $\cs$ is a suborder of
$Q_1\p Q_2$, in particular $P_i$ is also. By \ref{finest_chain_rep_Q},
$P_i$ is therefore a
suborder of $Q_1$ or $Q_2$. By \ref{interval_suborder},
if $P_i$ is a suborder of $Q_1$
then $P_1\p \cdots\p P_i$ is as well.
\ref{interval_suborder} similarly implies that if $P_i$ is a suborder of $Q_2$
then $P_i\p \cdots\p P_n$ is as also. The contrapositive is thus proved,
which completes the proof that $\forb(\cs)$ is
$\bit$-bit closed in this final case.

We now know that $\forb(\cs)$ is $S$-bit closed. Next, we show that every
$S$-bit closed set $X$ containing the empty and one point partial orders
has $\forb(P_1\p \cdots\p P_n)$ as a subset.

Suppose not. Then the $S$-bit closure $X$ of the set containing the empty
and one point partial orders is a proper subset of the $S$-bit closed
set $\forb(\cs)$.
Take a minimum cardinality SP order $Q$ in $\forb(\cs)$
that is not in $X$. Then $Q$ has at least two elements by choice of $X$.
Since $Q$ is an SP order, it follows that $Q$ is a chain or antichain sum.

If $Q$ is an antichain sum, then we may write $Q=Q_1\oplus Q_2$, where
$Q_1$ and $Q_2$ each have fewer elements than $Q$.  Since
$Q$ is a minimum size partial order in $\forb(\cs)-X$ by hypothesis,
we see that $Q_1$ and $Q_2$ are in $X$. Since $X$
is $\bit$-bit closed for $\bit$ the two point antichain $R_A$ with both points
labeled $R$, it follows that the antichain sum of two orders in $X$
is in $X$ as well. In particular, $Q$ is in $X$,
contrary to hypothesis. This contradiction shows that $Q$ can not be an
antichain sum.

Since $Q$ is not an antichain sum, $Q$ must be a chain sum $Q=Q_1\p Q_2$.
By choice of $Q$ as minimal, we know that $Q_1$ and $Q_2$ are in $X$.
Suppose $Q_2$ is in $\forb(P_n)$. Since $Q_1$ is in $X$ and $Q_2$ is
in $\forb(P_n)$, and since $X$ is $\bit$-bit closed for $\bit$
the two point chain with top labeled $\forb(P_n)$ and bottom labled $R$,
we see that $Q_1\p Q_2$ must be in $X$, contrary to hypothesis. Therefore
$Q_2$ is not in $\forb(P_n)$. By similar reasoning, $Q_1$ is not in
$\forb(P_1)$.

Choose the least $i$ such that $Q_1$ does not have a $P_1\p \cdots\p P_i$
suborder. Then $Q_1$ has a $P_1\p \cdots\p P_{i-1}$ suborder. If $Q_2$ has
a $P_i\p \cdots\p P_n$ suborder, then $Q_1\p Q_2$ has a $\cs$ suborder,
contrary to hypothesis. Therefore $Q_2$ has no $P_i\p \cdots\p P_n$ suborder.
Therefore $Q_1$ is in $\forb(P_1\p \cdots\p P_i)$ and
$Q_2$ is in $\forb(P_i\p \cdots\p P_n)$. Since the two point chain with
top labeled $\forb(P_i\p \cdots\p P_n)$ and bottom labeled $\forb(P_1\p \cdots\p P_i)$
is a bit in $S$ and $X$ is $S$-bit closed, it follows that $Q_1\p Q_2=Q$ is
in $X$, contrary to hypothesis.

In all cases, the assumption that $X$ is a proper subset of $\forb(\cs)$
is a contradiction. Equality therefore holds, thus completing the proof.
\end{prf}

To give a similar result for excluding a set of chain sums, we first need
some definitions.

\begin{dfn}
Fix $k\ge 1$.
For $1\le i\le k$ let $P_i$ be a chain sum.
A chain bit choice function for $(P_1,\ldots,P_k)$
is a function $c$ mapping each $P_i$ to a chain bit in
$\bss(P_i)$.
\end{dfn}

Given a chain bit $\bit$, we let $\b(\bit)$ and\\
$\t(\bit)$
denote the labels of the bottom and top points, respectively, of $\bit$.

In the next definition, we must intersect labels of bits. If all labels
are ideals, then no comment is necessary, but in general some labels may
be the symbol $R$, so we must extend the notion of intersection to include
this symbol. We make the convention that in the definition of bit set
corresponding to $(P_1,\ldots,P_k)$ below, the symbol $R$ is taken
to mean $\forb(P_1,\ldots,P_k)$. In other words, the intersection of
$R$ with a set is the intersection of $\forb(P_1,\ldots,P_k)$ and that set.
Moreover, if a rule tells us that a point should be labeled
$\forb(P_1,\ldots,P_k)$, we label that point $R$. Without this convention,
stating the following definition would be quite lengthy.

\begin{dfn}
Fix $k\ge 1$. For $1\le i\le k$, let $P_i$ be a chain sum.
The chain bit set $\bss(P_1,\ldots,P_k)$ corresponding to the tuple
$(P_1,\ldots,P_k)$ is the set of two point chain bits of the form
$$\bigcap_{1\le i\le k} \b(c(P_i))
\p
\bigcap_{1\le i\le k} \t(c(P_i)).$$
such that $c$ is a chain bit choice function for $(P_1,\ldots,P_k)$.
\end{dfn}

We note that the previous definition is consistent with 
\ref{chain_bit_choice_first_def} for the
case $k=1$. The following lemma generalizes \ref{one_chain_lem}
 to the case of excluding an arbitrary finite set of chain sums.

\begin{lem} \label{multi_chain_lem}
Let $k\ge 1$.
If the SP orders $P_1,\ldots,P_k$ are chain sums,
then $$\forb(P_1,\ldots,P_k)=L(\bss(P_1,\ldots,P_k)\cup\{R_A\}).$$
\end{lem}

\begin{prf}
For $k=1$, this is just \ref{one_chain_lem},
so we assume without loss of generality that $k\ge 2$.

Let $S=\bss(P_1,\ldots,P_k)\cup\{R_A\}$.
We must show that $\forb(P_1,\ldots,P_k)$ is the $S$-bit closure of the
doubleton containing the empty and one point partial orders.
Since $\forb(P_1,\ldots,P_k)$ trivially contains the empty and
one point partial orders, it is enough
to show that $\forb(P_1,\ldots,P_k)$ is $S$-bit closed and that every
$S$-bit closed set containing the empty and one point partial orders
has $\forb(P_1,\ldots,P_k)$ as a subset.

We first show that $\forb(P_1,\ldots,P_k)$ is $S$-bit closed.
By \ref{set_to_bit}, it is enough to show that $\forb(P_1,\ldots,P_k)$ is
$\bit$-bit closed for each bit $\bit$ in $S$.

First, if $\bit$ is $R_A$,
then to show that $\forb(\cs)$ is $\bit$-bit closed is simply to show that
$\forb(\cs)$ is closed under antichain sums. But this is exactly
\ref{antichain_bit_there}.

If $\bit\not=R_A$, then $\bit$ has the form
$$\bigcap_{1\le i\le k} \b(c(P_i))
\p
\bigcap_{1\le i\le k} \t(c(P_i))$$
for some chain bit choice function $c$ for $(P_1,\ldots,P_k)$.
To show that $\forb(P_1,\ldots,P_k)$ is $\bit$-bit closed is thus to
show that for each chain bit choice function $c$ for
$(P_1,\ldots,P_k)$, if $Q_1$ and $Q_2$ are SP orders
in $\forb(P_1,\ldots,P_k)$ such that $Q_1$ is in
$\bigcap_{1\le i\le k} \b(c(P_i))$
and $Q_2$ is in $\bigcap_{1\le i\le k} \t(c(P_i))$, then $Q_1\p Q_2$ is in
$\forb(P_1,\ldots,P_k)$ as well. To show that $Q_1\p Q_2$ is in
$\forb(P_1,\ldots,P_k)$, we must show that $Q_1\p Q_2$ forbids $P_i$
for $1\le i\le k$, so choose $i$.
Since $Q_1$ is in $\bigcap_{1\le i\le k} \b(c(P_i))$, in particular
$Q_1$ is in $\b(c(P_i))$. Similarly $Q_2$ is in $\t(c(P_i))$.
Since $c$ is a chain bit choice function for $(P_1,\ldots,P_k)$,
we see that $\b(c(P_i))\p \t(c(P_i))$ is a chain bit in $\bss(P_i)$.
Both $Q_1$ and $Q_2$ are in $\forb(P_i)$.
Therefore $Q_1\p Q_2$ is in $\forb(P_i)$ as needed. This completes
the proof that $\forb(P_1,\ldots,P_k)$ is $S$-bit closed.

We now know that $\forb(P_1,\ldots,P_k)$ is $S$-bit closed.
Next, we show that every
$S$-bit closed set $X$ containing the empty and one point partial orders
has $\forb(P_1,\ldots,P_k)$ as a subset.

Suppose not. Then the $S$-bit closure $X$ of the set containing the empty
and one point partial orders is a proper subset of the $S$-bit closed
set $\forb(P_1,\ldots,P_k)$.
Take a minimum cardinality SP order $Q$ in $\forb(P_1,\ldots,P_k)$
that is not in $X$. Then $Q$ has at least two elements by choice of $X$.
Since $Q$ is an SP order, it follows that $Q$ is a chain or antichain sum.

If $Q$ is an antichain sum, then we may write $Q=Q_1\oplus Q_2$, where
$Q_1$ and $Q_2$ each have fewer elements than $Q$.  Since
$Q$ is a minimum size partial order in $\forb(P_1,\ldots,P_k)-X$ by hypothesis,
we see that $Q_1$ and $Q_2$ are in $X$. Since $X$
is $\bit$-bit closed for $\bit$ the two point antichain $R_A$ with both points
labeled $R$, it follows that the antichain sum of two orders in $X$
is in $X$ as well. In particular, $Q$ is in $X$,
contrary to hypothesis. This contradiction shows that $Q$ can not be an
antichain sum.

Since $Q$ is not an antichain sum, $Q$ must be a chain sum $Q=Q_1\p Q_2$.
By choice of $Q$ as minimal, we know that $Q_1$ and $Q_2$ are in $X$.
For each $i$, since $Q_1\p Q_2$ is in
$\forb(P_i)=L(\bss(P_i)\cup\{R_A\})$, we know there is a two point chain
bit $B_i$ in $\bss(P_i)$ such that $Q_1$ is in
$\b(B_i)$ and $Q_2$ is in $\t(B_i)$. Define the chain bit choice function
$c$ for $(P_1,\ldots,P_k)$ by letting $c(P_i)=B_i$ for each $i$.
Then $Q_1$ is in $\bigcap_{1\le i\le k}\b(c(P_i))$ and $Q_2$ is in
$\bigcap_{1\le i\le k}\t(c(P_i))$. Moreover, $Q_1$ and $Q_2$ are in
$\forb(P_1,\ldots,P_k)$ and
$$\bigcap_{1\le i\le k}\b(c(P_i))\p
\bigcap_{1\le i\le k}\b(c(P_i))$$
is in $\bss(P_1,\ldots,P_k)$.
It follows that $Q=Q_1\p Q_2$ is in $\forb(P_1,\ldots,P_k)$,
contrary to hypothesis. This contradiction completes the proof.
\end{prf}

%something about formula and corresponding ideal

%We now define the formulas we use.
%We think of positive integers as {\em variables} of our formulas.
%The other symbols of our formulas are $\land$ and $\lor$, to be thought
%of as ``and'' and ``or'', respectively, and left and right parentheses.
%A finite sequence of symbols is called a string.
%Our formulas will be finite strings; that is, finite
%sequences of $($, $)$, $\land$, $\lor$, and positive integers.
%Given strings $f$ and $g$, we write $fg$ for the
%concatenation of $f$ and $g$, and we freely write expressions such
%as $(f\lor g)$ to mean the concatenation of string $($, $f$, $\lor$, $g$,
%and $)$ in that order.

%We define the set of formulas as the smallest such set of symbols
%satisfying the following conditions:
%\begin{enumerate}
%\item For all positive integers $i,j$, the strings
% $(i\land j)$ and $(i\lor j)$ are formulas.
%\item For all formulas $f$ and $g$, the strings $(f\land g)$
% and $(f\lor g)$ are formulas.
%\end{enumerate}

%Though $(1\land(2\land 3))$ and $((1\land 2)\land 3)$ are by definition
%distinct formulas, they are meant to and do represent the same thing.
%Moreover, though $1\land 2$ is by definition not a formula even though
%$(1\land 2)$ is, in practice, we only write the parantheses needed to
%eliminate ambiguity. We may write
%$1\land(2\lor 3)$ instead of $(1\land(2\lor 3))$, but we do not write
%$1\land 2\lor 3$ as this is ambiguous.
%With these conventions in place, no confusion should arise.

%We recursively define the ideal 

%\def\cf{\cal{F}}
%\def\cp{\cal{P}}
\def\otk{\{1,\ldots,k\}}
We now move onto excluding sets of antichain sums. As a motivating example,
we may wish to compute $\forb(P_1\oplus P_2,P_2\oplus P_3)$. We would then
let $\Gamma$ be the family of subsets of $\{1,2,3\}$
consisting of $\{1,2\}$ and $\{2,3\}$ and think of
$\forb(P_1\oplus P_2,P_2\oplus P_3)$ as
$$\bigcap_{F\in \Gamma}\forb\left(\bigoplus_{i\in F}P_i\right).$$

This example motivates us to define, given a sequence $P_1,\ldots,P_k$
of SP orders
%undone assume chain sum or singleton here?
and a family $\Gamma$ of nonempty subsets of $\otk$, the lower ideal
$$\forb(\Gamma;P_1,\ldots,P_k):=
\bigcap_{F\in \Gamma}\forb\left(\bigoplus_{i\in F}P_i\right)
.$$

We need several definitions.
A {\it splitting} of a set $X$ is an ordered pair
$(A,B)$ such that the sets $A$ and $B$ partition $X$.
We denote the set of splittings
of $X$ by $\s(X)$. A {\it  splitting function} for $X$ is a
function $h:\s(X)\to \{1,2\}$.

%Let $P_i$ be an SP order for each $i$ in a nonempty set $F$.
%Let $\ab$ be a splitting of $F$.
%If $Q_1$ is an SP order containing $\bo_{i\in A}P_i$
%and $Q_2$ is an SP order containing $\bo_{i\in B}P_i$, then
%$Q_1\o Q_2$  contains $\bo_{i\in F}P_i$.
%Therefore if $Q_1\o Q_2$ forbids $\bo_{i\in F}P_i$,
%then $Q_1$ forbids $\bo_{i\in A}P_i$ or
%$Q_2$ forbids $\bo_{i\in B}P_i$.
%This simple observation is at the heart of our structural description
%for $\fg$ below.

Let $\G$ be a family of subsets of $\otk$.
An {\it antichain bit choice function}, or ABCF, for $\G$
is a function $g$ with
domain $\G$ such that $g_F:=g(F)$ is a splitting function for $F$
for each set $F$ in $\G$.
We define the \lcis\ $\ln(g)$ of $g$ as the set of all pairs $(A,F)$ such that
$F$ is in $\G$ with $A\subseteq F$ and $g_F(A,F-A)=1$.
The \rcis\ $\rn(g)$ is defined similarly but with $g_F(A,F-A)=2$.

We define the left cell label $\lcl$ as the lower ideal
$$
\lcl:=
\fg
\cap
\bigcap_{(A,F)\in\ln(g)}\forb\left(\bo_{i\in A}P_i\right)
$$
and the right cell label $\rcl$ as the lower ideal
$$
\rcl:=
\fg
\cap
\bigcap_{(A,F)\in\rn(g)}\forb\left(\bo_{i\in F-A}P_i\right)
.$$

We now define $\bsspg$ as the set of labeled antichains that have the form
$$\lcl\o\rcl$$ for some ABCF $g$ for $\G$.

We need to use finest antichain partitions in the next lemma. This amounts
to assuming that our summands $\potk$ are not themselves antichain sums.

\begin{lem} \label{multi_antichain_lem}
Let $k\ge 1$.
If the SP orders $P_1,\ldots,P_k$ are not antichain sums,
then $$\fg=L(\bss(\G;P_1,\ldots,P_k)\cup\{R_C\}).$$
\end{lem}

\begin{prf}
Let $S=\bss(\G;P_1,\ldots,P_k)\cup\{R_C\}$.
We must show that $\fg$ is the $S$-bit closure of the
doubleton containing the empty and one point partial orders.
Since $\fg$ trivially contains the empty and
one point partial orders, it is enough
to show that $\fg$ is $S$-bit closed and that every
$S$-bit closed set containing the empty and one point partial orders
has $\fg$ as a subset.

We first show that $\fg$ is $S$-bit closed.
By \ref{set_to_bit}, it is enough to show that $\fg$ is
$\bit$-bit closed for each bit $\bit$ in $S$. 

First, if $\bit$ is $R_C$, then $\fg$ is $\bit$-bit closed by
\ref{chain_bit_there}.
Otherwise, by definition of $S$ and
$\bss(\G;P_1,\ldots,P_k)$, we see that $\bit$ must have the form
$\lcl\o\rcl$ for some ABCF $g$ for $\G$, so choose such a $g$.
To show that $\fg$ is $\bit$-bit closed for
$$\bit=\lcl\o\rcl,$$
we must show that if $Q_1$ is in $\lcl$ and $Q_2$ is in $\rcl$
then $Q_1\o Q_2$ is in $\fg$. Equivalently, we may show that
if $Q_1\o Q_2$ is not in $\fg$, then $Q_1$ is not in $\lcl$ or
$Q_2$ is not in $\rcl$.

Suppose $Q_1\o Q_2$ is not in
$$\forb(\Gamma;P_1,\ldots,P_k)=
\bigcap_{F\in \Gamma}\forb\left(\bigoplus_{i\in F}P_i\right)
.$$
Then there is $F$ in $\G$ such that $Q_1\o Q_2$ is not in
$\forb\left(\bigoplus_{i\in F}P_i\right)$. Therefore
$Q_1\o Q_2$ contains a $\bigoplus_{i\in F}P_i$ suborder.
We may then choose a one to one order preserving map
$h:\bo_{i\in F}P_i\to Q_1\o Q_2$ embedding $\bo_{i\in F}P_i$
into $Q_1\o Q_2$.
Since no $P_i$ is an antichain sum, we know by
\ref{antichain_bit_there}
that $h(P_i)$ is contained in $Q_1$ or $Q_2$ for each $i$.
Let $A=\{i\in F: h(P_i)\subseteq Q_1\}$.
Then $F-A=\{i\in F: h(P_i)\subseteq Q_2\}$.
If $A$ is empty then $\bo_{i\in F}P_i$ is a suborder of $Q_2$.
Therefore $Q_2$ is not in
$\forb\left(\bo_{i\in F}P_i\right)$, which implies $Q_2$ is not
in
$$\bigcap_{F\in \G}\forb\left(\bo_{i\in F}P_i\right).$$
By the definition of $\rcl$, this in turn implies that
$Q_2$ is not in $\rcl$.
This proves our claim in the case that
$A$ is empty. Similarly if $F-A$ is empty. We may thus assume
that $A$ and $F-A$ are nonempty.

Either $g_F(A,F-A)=1$ or $g_F(A,F-A)=2$. If $g_F(A,F-A)=1$, then
$(A,F)$ is in $\ln(g)$. Certainly $\bo_{i\in A}P_i$ is not in
$\forb\left(\bo_{i\in A}P_i\right)$, and $Q_1$ contains $\bo_{i\in A}P_i$,
which implies $Q_1$ is not in $\forb\left(\bo_{i\in A}P_i\right)$.
Therefore $Q_1$ is not in
$$\bigcap_{(A,F)\in\ln(g)}\forb\left(\bo_{i\in A}P_i\right).$$ By definition
of $\lcl$, we thus see that $Q_1$ is not in $\lcl$.
Similarly, if $g_F(A,F-A)=2$ then $Q_2$ is not in $\rcl$, as was to be
shown. This completes the proof of the claim that $\fg$ is $S$-bit closed.

We must now show that every $S$-bit closed set containing the empty
and one point partial orders has $\fg$ as a subset. Suppose not.
Then the $S$-bit closure $X$ of the set containing the empty
and one point partial orders is a proper subset of the $S$-bit closed
set $\fg$.
So take a minimum cardinality SP order $Q$ in $\fg$
that is not in $X$. Then $Q$ has at least two elements by choice of $X$.
Since $Q$ is an SP order, it follows that $Q$ is a chain or antichain sum.
If $Q$ is a chain sum $Q_1\p Q_2$ then $Q_1$ and $Q_2$ are in $X$ by
choice of $Q$ as minimal. Since $R_C$ is in $S$ and $X$ is $S$-bit
closed, it then follows that $Q=Q_1\p Q_2$ is in $X$, contrary to
hypothesis. This contradiction shows that $Q$ is an antichain sum.

We write $Q=Q_1\o Q_2$.
We wish to get a contradiction in this case as well by showing in fact
that $Q$ is in $X$. Since $Q_1$ and $Q_2$ are in $X$ by minimality of $Q$,
and since $X$ is $\bit$-bit closed for
$$\bit=\lcl\o\lcl,$$
we see it is
enough to show there is an ABCF $g$ for $\G$ such that $Q_1$ is in $\lcl$
and $Q_2$ is in $\rcl$. Since $Q$ is in the lower ideal $\fg$, the
suborders $Q_1$ and $Q_2$ are in $\fg$ as well. By definition of $\lcl$
and $\rcl$, it is therefore enough to exhibit an ABCF $g$ for $\G$ such
that $Q_1$ is in
$$\bigcap_{(A,F)\in\ln(g)}\forb\left(\bo_{i\in A}P_i\right)$$
and
$Q_2$ is in
$$\bigcap_{(A,F)\in\ln(g)}\forb\left(\bo_{i\in F-A}P_i\right).$$

Choose $F$ in $\G$. Since $Q_1\o Q_2$ is in $\fg$, we see that $Q_1\o Q_2$
forbids $\bo_{i\in F}P_i$. Therefore for each splitting $(A,B)$ of $F$,
the SP order $Q_1$ must forbid  $\bo_{i\in A}P_i$ or $Q_2$ must forbid
$\bo_{i\in B}P_i$. Consider the ABCF $g$ for $\G$ such that for each
$F$ in $\G$ and each splitting $(A,B)$ of $F$, we have
$g_F(A,B)=1$ if $Q_1$ forbids $\bo_{i\in A}P_i$ and
$g_F(A,B)=2$ otherwise.

To show that $Q_1$ is in
$$\bigcap_{(A,F)\in\ln(g)}\forb\left(\bo_{i\in A}P_i\right),$$
it is enough to show that $Q_1$ is in $\forb\left(\bo_{i\in A}P_i\right)$
for each $F$ in $\G$ and each nonempty $A\subseteq F$ such that
$g_F(A,F-A)=1$. This is immediate from the definition of $g_F$.
Similarly, it follows immediately from the definition of $g_F$ that
$Q_2$ is in
$$\bigcap_{(A,F)\in\ln(g)}\forb\left(\bo_{i\in F-A}P_i\right).$$
This completes the proof of the lemma.
\end{prf}

%undone is \bss(A) for set A defined yet?
%is it convention that everything is SP order? should be
\begin{lem} \label{mixed_case_lem}
If $A$ and $B$ are nonempty sets of chain sums and antichain sums,
respectively, then $\forb(A\cup B)=L(\bss(A)\cup\bss(B))$.
\end{lem}

\begin{prf}
We know that
$\forb(A)$ is $\bit$-bit closed for each bit $\bit$ in $\bss(A)$. We also
know by \ref{antichain_bit_there}
that $\forb(A)$ is closed under arbitrary antichain sums,
and since each bit in $\bss(B)$ is an antichain, we see
that $\forb(A)$ is $\bit$-bit closed for each $\bit$ bit in $\bss(B)$.
Therefore $\forb(A)$ is $\bit$-bit closed for each bit $\bit$ in
$\bss(A)\cup \bss(B)$. By similar reasoning, $\forb(B)$ is $\bit$-bit closed
for each bit $\bit$ in $\bss(A)\cup \bss(B)$ as well.
This implies that $\forb(A\cup B)=\forb(A)\cap\forb(B)$ is $\bit$-bit
closed for each bit $\bit$ in $\bss(A)\cup \bss(B)$, and hence
$\forb(A\cup B)$ is $\bss(A)\cup \bss(B)$ closed. Therefore
$L(\bss(A)\cup\bss(B))\subseteq \forb(A\cup B)$.

If $\forb(A\cup B)=L(\bss(A)\cup\bss(B))$, we are done. Suppose not.
Then $L(\bss(A)\cup \bss(B))$ is a proper subset of  $\forb(A\cup B)$ .
Choose a minimum cardinality SP order $Q$ in $\forb(A\cup B)$
that is not in  $L(\bss(A)\cup \bss(B))$. Since $Q$ has at least two points,
$Q$ is a chain sum or an antichain sum. We assume that $Q$ is a chain
sum. The case that $Q$ is an antichain sum is entirely similar.

Since $Q\in\forb(A\cup B)\subseteq\forb(A)$, we see that $Q$ is in
$\forb(A)=L(\bss(A)\cup\{R_A\})$. Therefore there is a bit $\bit$ in
$\bss(A)\cup\{R_A\}$ that generates $Q$ from proper suborders.
Since $Q$ is a chain sum, we know that $Q$ is not an antichain sum.
Therefore $\bit\not=R_A$, which implies $\bit$ is in $\bss(A)$. In
particular, the $\bss(A)\cup\bss(B)$-bit closure of the set of
proper suborders of $Q$ contains $Q$. Since each proper suborder of $Q$
is in $L(\bss(A)\cup\bss(B))$ and $L(\bss(A)\cup\bss(B))$ is
$\bss(A)\cup\bss(B)$-bit closed, we see that $Q$ is in
$L(\bss(A)\cup\bss(B))$, contrary to assumption. This contradiction
completes the proof.
\end{prf}

\section{The Main Theorem}

\begin{thm}
There is a structural description for each nontrivial proper lower ideal $L$.
\end{thm}

\begin{prf}
The proper lower ideal $L$ is described by a finite list of forbidden
suborders.
That list either consists of one chain sum, multiple chain sums, multiple
antichain sums, or both chain and antichain sums.
We thus use \ref{one_chain_lem}, \ref{multi_chain_lem},
\ref{multi_antichain_lem}, or \ref{mixed_case_lem}, respectively to obtain a
set $S$ of bits generating $L$. Each label of a partial order in $S$ is
either the symbol $R$ or is an ideal properly contained in $L$. For properly
contained ideals, we repeat this procedure recursively. We thus obtain a
finitely branching tree representing this construction. By the fact that
SP orders are better quasi ordered under the
suborder relation \cite{countable_sp_bqo}, it follows
that there is no infinite descending sequence of lower ideals of SP orders.
Therefore this construction tree is a finitely branching tree with no
infinite branch, which is finite by K\"onig's Lemma.
This completes the proof.
\end{prf}

We stress that this theorem is not just theoretical; it can be applied by
hand in practice to obtain specific structure theorems quickly. As one example,
we characterize the diamond free SP orders.
The diamond is the unique poset on points $a,b,c,d$ such that $a<b<d$,
$a<c<d$, and $b$ and $c$ are incomparable.
An SP order is called diamond free if there is no diamond suborder.
A (partial order theoretic) tree is a poset such that for each $x$, there
are no incomparable elements less than $x$. A forest is tree or an antichain
sum of trees. An upside down tree (forest) is a poset such that the reverse
order is a tree (forest).
A forest on top of an upside down forest is a chain sum of a forest and
upside down forest with the outer poset a two point chain, the top poset
a forest, and the bottom poset an upside down forest.
With these definitions, the reader may use the
results of this paper to quickly prove the following corollary. 

\begin{cor}
A finite SP order is diamond free iff it has the form
$$\lexpi,$$
where $(I,\le_I)$ is an antichain and $P_i$ is a forest on top of
an upside down forest for each $i$.
\end{cor}

Note that the structural descriptions for ideals are
not at all in general unique. Our procedure simply finds one of them.
The one found may in fact have redundant rules.
Note also that since lemmas \ref{one_chain_lem}, \ref{multi_chain_lem},
\ref{multi_antichain_lem}, and \ref{mixed_case_lem}, only involve
the two point chain and antichain $R_A$ and $R_C$, it follows that each lower
ideal has a structural description only involving two point posets at any
depth. At least to the author, this fact was initially surprising.

\section{Acknowledgements}

I thank Yared Nigussie for teaching me the mathematics
\cite{nigussie}.
The deeper understanding thus obtained by the author made this paper
possible.

\bibliographystyle{plain}
\bibliography{mybib}
\end{document}